\documentclass{article}
\usepackage{fullpage}
\usepackage[utf8]{inputenc} 
\usepackage[T1]{fontenc}    
\usepackage{url}            
\usepackage{booktabs}       
\usepackage{amsfonts}       
\usepackage{nicefrac}       
\usepackage{microtype}      
\usepackage{lipsum}
\usepackage{graphicx}       
\graphicspath{{media/}}     
\usepackage{amsmath, amssymb}
\usepackage{longtable}
\usepackage{makecell}
\usepackage{xcolor}
\usepackage{float}
\usepackage{enumitem}
\usepackage[bookmarks=false]{hyperref}
\usepackage{pstricks}
\def\Y{Y}
\def\x{\mathbf x}
\def\u{\mathbf u}
\def\a{\mathbf a}
\def\b{\mathbf b}
\def\xx{\bar{\x}}
\def\uu{\bar{\u}}
\def\c{\mathbf c}
\def\cc{\bar{\c} }
\def\thc{\theta}
\def\varthc{\vartheta}
\def\qc{q}
\title{Symmetry analysis and new partially invariant solutions for the gas dynamics system with a special equation of state
}

\author{
	Dilara Siraeva  and Irina Kogan\\
	Department of Mathematics\\ North Carolina State University, Raleigh\\
	\texttt{sirdilara@gmail.com and iakogan@ncsu.edu} 
}

\begin{document}
\maketitle
\begin{abstract}
This paper is a contribution to the symmetry analysis of the gas dynamics system in the vein of the ``podmodeli'' (submodels) program outlined by Ovsyannikov (1994). We consider the case of the special state equation,  prescribing pressure to be the sum of entropy and an arbitrary function of density.  Such a system has a  12-dimensional symmetry Lie algebra. This work advances the study of its four-dimensional subalgebras, continuing the work started in Siraeva (2024). For a large subset of not previously considered, non-similar four-dimensional subalgebras from an optimal list in Siraeva (2014), we compute a complete set of generating invariants. For one of the subalgebras, we construct a partially symmetry-reduced system. We explicitly solve this reduced system (submodel). This leads to new families of explicit solutions of the original system. We analyze the trajectories of these solutions.  Additionally, we match each of the subalgebras considered in this paper with its isomorphism class, planting a seed for future study of the hierarchy of the reduced systems.

\end{abstract}

\noindent {\bf Keywords}: Gas dynamics system, special state equation, symmetry Lie algebras, classification, invariants, symmetry-reduction, submodel, exact solutions. 

\noindent {\bf MSC 2020:} 35B06, 35A30.
\section{Introduction}
This paper is a contribution to an extensive and ongoing research program of the symmetry analysis   of the system of  gas dynamics equations:
\begin{align}\label{eq-gd}
	\nonumber D\,\u+\rho^{-1}{\bf\nabla} P&=0,\\
        D\,\rho+\rho\,{\rm div}\,\u&=0,\\
        \nonumber D\,P+\rho f_{\rho}{\rm div}\,\u&=0,
\end{align}
where $D=\partial_t+(\u\cdot{\bf\nabla})$ is a total differentiation operator, $t$ is time; $\nabla=\partial_{\x}$ is the gradient with respect to the position vector $\x\in\mathbb R^3$;   $\u\in \mathbb R^3$ is the velocity vector; $\rho$ is density; $P$ is  pressure, and $f(\rho,S)$ is an arbitrary function of density and entropy, prescribing the dependency of $P$ on these quantities via the \emph{state equation:}  
\begin{equation}\label{state-arb}
	\begin{array}{c}
P=f(\rho,S).
\end{array}
\end{equation}
The subscript denotes the derivative with respect to the corresponding variable.
Under the assumption that {$f_S\neq 0$}, the last equation of system \eqref{eq-gd} can be replaced by
\[DS=0.\]
Due to its prominent role in continuous mechanics, finding new exact solutions to system \eqref{eq-gd}-\eqref{state-arb} is an important problem.
One of the most effective tools for finding such solutions is the symmetry analysis method pioneered by Sophus Lie \cite{Lie1881, Lie1883, Lie1885} and further developed in \cite{Ovs1982, Olver1986, Bluman1989, Ibrag1994, Ibrag2001, Hydon2000, Arrigo2015} and other works. This method uses symmetries and invariants to reduce the original system to an easier one to solve. Each solution of the reduced system can be lifted to a family of solutions of the original system.   

  In this paper,  we  closely follow  the program outlined by Lev Ovsyannikov \cite{Ovs1994}. As shown  in \cite{Ovs1982}, the full symmetry group of \eqref{eq-gd}-\eqref{state-arb} with arbitrary $f$ is an 11-dimensional Lie group generated by   continuous transformations: time translations, space translations, Galilean translations, rotations and uniform dilatation of space-time; and discrete transformations: change orientation and  time inversions, as summarized in Table~\ref{table_sym} of this paper. For a specific $f$, the symmetry group can be larger, as shown in Table 1 of \cite{Ovs1994}. In this paper, we consider, one of the less studied  cases, where   
\begin{equation}\label{state-sp}
	\begin{array}{c}
		P=f(\rho)+S,
	\end{array}
\end{equation}
admitting an additional one parametric group of symmetries of pressure translations. Although the authors are not aware of any physical phenomena modeled by this state equation, we think that completing Ovsyannikov's program is a worthwhile mathematical goal, and we hope that (as it happened before with many other "purely mathematical" constructions) some venues of application may present themselves in the future.

 Ovsyannikov's program consists of finding a complete list of the subgroups of the full symmetry group, classified up to a conjugation by an element of the full group (\emph{inner automorphisms}) and, in some cases, up to some additional \emph{outer automorphisms}. Two subgroups in the same equivalence class are called \emph{similar}. It is customary, in the area of symmetry reduction, to work infinitesimally, replacing a symmetry group with the corresponding symmetry Lie algebra.
A basis of the 12-dimensional Lie algebra, denoted $L_{12}$, of the symmetry group of   system~\eqref{eq-gd} with the special state equation \eqref{state-sp}  is given in Section~\ref{sect-alg}.
 It is well known that invariant functions are annihilated by the infinitesimal generators.

 Similarity of subgroups of the full symmetry group translates into similarity of Lie subalgebras of the full symmetry Lie algebra.  A list consisting of a representative  of each similarity class is called an \emph{optimal list}. The Ovsyannikov's program entails finding such an optimal list and then computing a  generating set of independent invariant functions  for each subalgebra in this list. In  this paper, we consider \emph{regular reductions}, in the terminology of \cite{Khab2019,Ovs1995}, where we maximize the number of  generating invariants that depend only on the ``old'' independent variables $(t,x,y,z)$, and assign all such invariants to be ``new'' independent variables. The rest are  chosen to be new dependent variables.

An  \emph{invariant submodel} is a reduced system, obtained by rewriting the original system in terms of these new  invariant dependent and independent variables. \emph{A submodel is partially invariant} if we include one or more  non-invariant functions among the list of new independent variables. The number of such additional non-invariant independent variables is called the \emph{defect of the submodel}.  The number of new independent variables is called the \emph{rank of the submodel}.   Any explicit solution for a submodel leads to a family of explicit solutions for the original system.

In \cite{Ovs1994}, Ovsyannikov admits that his program is ``completely realistic but extremely laborious, and its use requires a great deal of collective work.''    Although computation can be facilitated by software, it requires significant ``by-hand'' involvement, laborious case-by-case analysis, which involves some dose of tricks and ingenuity. 
The current paper contributes the following {\bf new results} for system \eqref{eq-gd} with the state equation given by \eqref{state-sp}: 
\begin{itemize}
\item In Section~\ref{sect-inv}, we present generating sets of invariants for a large number of non-similar four-dimensional subalgebras of the full symmetry Lie algebra $L_{12}$ of \eqref{eq-gd}, \eqref{state-sp}.   The subalgebras are numbered according to the optimal system presented in \cite{Sir2014}. Some of the numbered items in the list correspond to parametric families of non-similar subalgebras. In certain cases, it turned out to be more efficient to compute and present invariants in cylindrical or spherical coordinates, rather than in the standard Cartesian coordinates. The invariants were computed and verified with the help of {\sc maple} software. 
\item In Section~\ref{sect-isom}, we matched each four-dimensional subalgebra considered in Section~\ref{sect-inv} with its isomorphism type according to the classification of 4-dimensional Lie algebras presented in \cite{PatWint1976, PatWint1977}. (It is important to keep in mind that non-similar Lie algebras can be isomorphic.) This result, combined with the classification of subalgebras of 4-dimensional Lie algebras \cite{PatWint1977}, can significantly facilitate an important task of obtaining a hierarchy of the submodels.
\item Each of the generating sets of invariants, computed in \ref{sect-inv}, can be used to obtain a submodel. 
In Section~\ref{sect-sol}, we derive a submodel, of rank one and defect one, corresponding to the last subalgebra (labeled 4.77) in our list. We then find explicit solutions for this submodel, which, in turn, lead to new solutions of the gas dynamics system \eqref{eq-gd}, \eqref{state-sp}. We analyzed the trajectories of these new solutions. 

\end{itemize}

We now put our results in the  context of {\bf previous works}. 
An optimal list  of subalgebras of the  symmetry Lie algebra $L_{11}$ of system \eqref{eq-gd} with an arbitrary state equation \eqref{state-arb} appeared in \cite{Ovs1994}. The graph of embedded subalgebras was constructed in \cite{MukmKhab2019}, and there is a vast body of works on computing  the corresponding invariants, submodels, and solutions \cite{ChirKhab2012,Khab2013,Mamont1999,Mamont2007,Ovs1996,Ovs1999-results}.  In particular, four-dimensional subalgebras of $L_{11}$ were extensively studied  in~\cite{Khab2019}, where generating sets for all 48 types of the four-dimensional subalgebras were computed  and some new regular and irregular partially invariant solutions were obtained.

The specialized equations of state for which an extension of the 11--dimensional Lie algebra occurs were listed in \cite[Chapter~3]{Ovs1982} and \cite{Ovs1994}. In~\cite{Khab2011}, these Lie algebras were classified up to isomorphism.  For each Lie algebra from this list,  optimal lists of subalgebras were computed in~\cite{Cherevko1996, Golovin1996, Khab2014, Makarevich2011, Ovs1994, Sir2014}. For some of the subalgebras, the symmetry reduction was carried out completely, leading to new explicit  solutions. In particular,    
for gas dynamics equations with a polytropic gas state equation $P=S
\rho ^{\gamma}$, $\gamma>0$, simple solutions, i.e., solutions of rank and defect equal to zero,
were  obtained in~\cite{Ovs1999-simple}. For the state equation of a monatomic gas {$P=S
\rho ^{5/3}$}, eight simple invariant solutions  were obtained from four--dimensional subalgebras were constructed in \cite{Nikon2023}.

In the current paper, the arbitrary state equation \eqref{state-arb} is specialized to \eqref{state-sp}, and so the symmetry group is enlarged to $L_{12}$. The optimal list of subalgebras (excluding, of course,  subalgebras of $L_{11}$ and trivial direct products with a subalgebra from $L_{11}$) was first found in \cite{Sir2014}. The current paper is a continuation of the study, initiated in  \cite{Sir2024-partially}, of \emph{four-dimensional} subalgebras from this list. In \cite{Sir2024-partially}, invariants for 29  items (some including parameters)  were calculated.   With 20 more items considered here, there are still 28 four-dimensional subalgebras whose invariants remain to be computed in future works.

As far as the lower-dimensional  Lie subalgebras of $L_{12}$ are concerned, the following results were previously obtained. Among two-dimensional subalgebras, only two define partially invariant submodels of rank 3 and defect 1. Their reductions to invariant submodels were obtained in~\cite{Sir2018-reduction}. The rest of the two-dimensional subalgebras admit invariant submodels of rank 2. In ~\cite{Sir2019-canonic, Sir2019-classif}, all these submodels were explicitly written in the canonical form (see \cite{Ovs1999-results} for the definition of canonical forms). In ~\cite{SirKhab2018}, explicit solutions for one of such submodels were computed and analyzed.
For the three-dimensional subalgebras of $L_{12}$, exact solutions were obtained in \cite{NikSirYulm2022, Sir2020, Sir2021} for four submodels of rank 1.   Using three-dimensional subalgebra generated by space translations, Galilean translations, and pressure translation, a family of exact solutions for the system was obtained in ~\cite{Sir2024-invariant}, which describes the motion of particles with a linear velocity field and non-homogeneous deformation in the 3D-space.

\section{The symmetry Lie algebra and its automorphisms}\label{sect-alg}

Equations  \eqref{eq-gd} represent the gas-dynamic system in the Cartesian (Descartes') coordinate system with four independent variables 
\begin{equation}\label{d-i}
   t \text{ and } \x=(x,y,z)
\end{equation}
and five dependent variables
\begin{equation}\label{d-d}
\u=(u,v,w),\quad \rho, \text{ and } P.
\end{equation}
The symmetry group of  \eqref{eq-gd} is the largest group of diffeomorphisms on  the 9-dimensional space, parameterized by $t$, coordinates of the position and velocity vectors, pressure and density, that maps each solution of \eqref{eq-gd}  to another solution. 
With an arbitrary equation of state~\eqref{state-arb}, the full symmetry group is generated  by $11$ one-parametric group of transformations and  two involutions listed in Table~\ref{table_sym}.  The corresponding basis of  the 11-dimensional symmetry Lie algebra  $L_{11}$, relative to the Cartesian coordinates, $\x=(x,y,z)$ and $\u=(u,v,w)$ is listed in the right column of the table. Here, as usual $SO_3$ denotes the special orthogonal group of rotations in $\mathbb R^3$ and $\mathbb R^*$ denotes non-zero real numbers.

\begin{table}[!h]
  \begin{center}
\caption{Symmetries of \eqref{eq-gd}-\eqref{state-arb}}\label{table_sym}
\begin{tabular}{|l|l|l|}
\hline
&Transformations& Infinitesimal Generators\\
\hline
 Space translations (ST) &$\xx=\x+\a,\qquad \a\in \mathbb R^3$ & $X_1=\partial_{x}, \quad X_2=\partial_{y}, \quad X_3=\partial_{z} $\\ 
\hline
Galilean translations (GT) & $\xx=\x+t\,\b,\quad \uu=\u+\b $
& $ X_4=t\partial_{x}+\partial_{u},\quad X_5=t\partial_{y}+\partial_{v}$, \\
&  $\b\in \mathbb R^3$ & $X_6=t\partial_{z}+\partial_{w}$ \\ 
\hline
Rotations (R) &	$\xx=R\x,\quad \uu=R\u, $ & $X_7=y\partial_{z}-z\partial_{y}+v\partial_{w}-w\partial_{v}$,\\
& & $X_8=z\partial_{x}-x\partial_{z}+w\partial_{u}-u\partial_{w}$,\\
&$R\in SO_3$ &$X_9=x\partial_{y}-y\partial_{x}+u\partial_{v}-v\partial_{u}$\\
\hline
Time translation (TT) &	$\bar t=t+\tau,\quad \tau\in\mathbb R $  & $X_{10}=\partial_{t} $ \\ 
\hline
Uniform dilation (D) & $\bar t=\lambda\,t,\quad \xx=\lambda\,\x,\quad \lambda \in\mathbb R^*$ & $X_{11}=t\partial_{t}+x\partial_{x}+y\partial_{y}+z\partial_{z}$, \\
\hline
Change of the orientation ($I_1$) & $\x^{\,*}=-\x,\quad \u^{\,*}=-\u$ &\\
\hline
Time inversion ($I_2$) & $\bar t=-{t},\quad \uu=-\u$ &\\
\hline
\end{tabular} 
\end{center}
\end{table}

The symmetry group of equations~~\eqref{eq-gd} with special equation of state \eqref{state-sp} is generated by an additional one-parametric symmetry group: the pressure translation~\cite{ChirKhab2012}:
\begin{table}[!h]
  \begin{center}
\caption{Additional symmetries of \eqref{eq-gd} with the special state equation \eqref{state-sp}}\label{table_sym_ext}
\begin{tabular}{|l|l|l|}
\hline
&Transformations& Infinitesimal Generator\\
\hline
Pressure translation &$P^{*}=P+c,\qquad c\in \mathbb R$ & $\quad \Y=\partial _P. $\\ 
\hline
\end{tabular} 
\end{center}
\end{table}

The commutators of  the basis of infinitesimal generators $X_1,\dots X_{11}$ of $L_{11}$, listed in Table~\ref{table_sym}, are given in 
Table~\ref{tabl_komm} \cite{Ovs1994}, where instead of generators $X_i$, $i=1,\dots,11$, we simply write indices $i$. Empty entries represent zeros. The additional generator $\Y=\partial _P$ commutes with  $X_i$, $i=1...11$ and thus Table~\ref{tabl_komm}  contains all non-trivial commutators relationships for $L_{12}$.
\begin{table}[!h]
	\caption{\mbox{The table of commutators of infinitesimal generators of Lie algebra~$L_{11}$}.}\label{tabl_komm}
	\begin{center}
		\begin{tabular}{ c p{2.4cm} p{2.4cm} p{2.4cm} p{1.7cm}  }
			\hline
			{} & $\hspace{0.3cm}{1}$ $\hspace{0.5cm}{2}$ $\hspace{0.5cm}{3}$     & ${\hspace{0.3cm}4}$ ${\hspace{0.5cm}5}$ ${\hspace{0.5cm}6}$   &  $\hspace{0.3cm}{7}$  $\hspace{0.5cm}{8}$ $\hspace{0.5cm}{9}$ & $\hspace{0.2cm}{10}$ $\hspace{0.4cm}{11}$ \\   \hline
			
			{1} &                 &                 & $\hspace{0.85cm}{-3}$  $\hspace{0.4cm}{2}$     &  $\hspace{1.3cm}{1}$  \\
			
			{2} &                 &                 & $\hspace{0.3cm}{3}$    $\hspace{1cm}{-1}$    &  $\hspace{1.3cm}{2}$  \\
			
			{3} &                 &                 & $\hspace{0.03cm}{-2}$   $\hspace{0.4cm}{1}$     &  $\hspace{1.3cm}{3}$  \\  
			
			{4} &                 &                 & $\hspace{0.85cm}{-6}$  $\hspace{0.4cm}{5}$     &  $\hspace{0.09cm}{-1}$  \\
			
			{5} &                 &                 & $\hspace{0.3cm}{6}$    $\hspace{0.95cm}{-4}$    &  $\hspace{0.09cm}{-2}$  \\
			
			{6} &                 &                 & $\hspace{0.05cm}{-5}$   $\hspace{0.4cm}{4}$      & $\hspace{0.09cm}{-3}$                 \\   
			
			{7} &  $\hspace{0.85cm}{-3}$ $\hspace{0.4cm}{2}$ & $\hspace{0.78cm}{-6}$ $\hspace{0.47cm}{5}$ & $\hspace{0.85cm}{-9}$ $\hspace{0.4cm}{8}$   &    \\
			
			{8} &  $\hspace{0.3cm}{3}$ $\hspace{0.95cm}{-1}$  & $\hspace{0.3cm}{6}$ $\hspace{0.95cm}{-4}$ &  $\hspace{0.3cm}{9}$ $\hspace{0.95cm}{-7}$   &    \\
			
			{9} &  $\hspace{0.03cm}{-2}$ $\hspace{0.5cm}{1}$  &    $\hspace{0.03cm}{-5}$ $\hspace{0.45cm}{4}$ & $\hspace{0.045cm}{-8}$ $\hspace{0.5cm}{7}$ & \\ 
			
			{10}&                  &   $\hspace{0.3cm}{1}$  $\hspace{0.47cm}{2}$ $\hspace{0.47cm}{3}$ &           & $\hspace{1.1cm}{10}$  \\
			{11}& $\hspace{0.03cm}{-1}$ $\hspace{0.09cm}{-2}$ $\hspace{0.09cm}{-3}$ &  &  &  $\hspace{0cm}{-10}$  \\   \hline
			
		\end{tabular}
	\end{center}
\end{table}

Each of the transformations in Table~\ref{table_sym} gives rise to an inner automorphisms of $L_{12}$.  An arbitrary $X\in L_{12}$ can be written as
$$X=c_0\,Y+\sum_{i=1}^{11} c_i\, X_{i},$$
for arbitrary $c_i\in\mathbb R$, $i=0,\dots,11$.
The image of $X$ under an automorphism, is then
$$\bar X=\bar c_0\,Y+\sum_{i=1}^{11} \bar c_i\, X_{i},$$
where  formulas for  $\bar c_i$, $i=0,\dots,11$ in terms of $c_i$, $i=0,\dots,11$ and transformation parameters is given in Table~\ref{table_autom}~\cite{Ovs1994}. To shorten the formulas, we  group $c_1,\dots, c_9$ (and their transformed versions) into three tuples of three:  
$$\c_1=(c_1,c_2,c_3),\quad \c_2=(c_4,c_5,c_6), \quad \c_3=(c_7,c_8,c_9)$$
Coefficients, whose transformations do not appear in the table are invariant (unchanged). 
\begin{table}[!h]
  \begin{center}
    \caption{Inner Automorphisms of Lie algebra $L_{12}$}\label{table_autom}
	\begin{tabular}{|c|p{9.5cm}|}
	\hline
	ST & $\cc_{1}=\c_1+c_{11}\,\a-\a\,\times \c_{3},$\hfill $\quad \a\in\mathbb R^3$\\      \hline
	GT & $\cc_1=\c_1-c_{10}\,\b,\quad  \cc_2=\c_2-\b\,\times \c_3$,\hfill $\quad             \b\in\mathbb R^3$\\ 
        \hline
	R & $\cc_{1}=R\,\c_1,\quad\cc_{2}=R\,\c_2,\quad \cc_{3}=R\,\c_3,$ \hfill $\quad       R\in SO_3$ \\ 
        \hline
			$TT$ &  $\cc_{1}=\c_{1}+\tau\, \c_{2},\quad \bar c_{10}=c_{10}+\tau\, c_{11}$,  \hfill $\quad\tau\in\mathbb R$ \\ \hline
			D & $\cc_{1}=\lambda\, \c_{1},\quad  \bar c_{10}=\lambda\, c_{10}$ \hfill $\quad\lambda\in\mathbb R^*$ \\ \hline
			$I_{1}$ & $\cc_{1}=-\c_{1},\ \quad \cc_{2}=-\c_{2}$ \\ 
            \hline
			$I_{2}$ & $\cc_{2}=-\c_{2},\ \quad \bar c_{10}=-c_{10}$ \\ \hline
		\end{tabular} 
  \end{center}
		\end{table}
Since pressure translations commute with all transformations in Table~\ref{table_sym}, they do not induce any inner automorphisms of $L_{12}$. However,  it is easy to check that scaling of the coefficient in front of $Y$:
\begin{equation}\label{outer}\bar c_0=\mu\,  \c_0, \quad \mu\in \mathbb R^*
\end{equation}
is an outer  automorphisms of $L_{12}$.

Two subalgebras of $L_{12}$ are \emph{similar} if  one can be mapped to another by a composition of inner from   outer automorphism given in Table~\ref{table_autom}  and formula \eqref{outer} respectively. An optimal list consists of one representative from each class. Since invariants of similar subalgebras are easily related to each other, one is interested only in computing invariants from the optimal list. An optimal list for subalgebras of  $L_{12}$ was first computed in \cite{Sir2014}.

An optimal list of  for $L_{12}$ subalgebras was computed in \cite{Sir2014} using the following considerations. Let $L\in L_{12}$ be a Lie subalgebra of dimension $2\leq n\leq 11$. Since  $L_{12}=L_{11}\oplus\{Y\}$, it is not difficult to show that $\hat L =L\cap L_{11}$ has either dimension $n$  or $n-1$. In the first case, $L=\hat L\subset L_{11}$ is similar to a subalgebra  from the  optimal list for  subalgebras of $L_{11}$ appearing in \cite{Khab2013}.  In the second case $L$ has a basis $Z_1,\dots ,Z_{n-1}, Z_n+Y$, where $Z_1,\dots ,Z_{n-1}$ is a basis of an $(n-1)$-dimensional subalgebra of $\hat L\in L_{11}$. Up to inner automorphisms listed in Table~\ref{table_autom}, we can assume  that it is one of the $L_{11}$-subalgebras appearing in the optimal list from \cite{Khab2013}. Therefore, to find an optimal list of $n$-dimensional subalgebras of $L_{12}$,  one can start by appending a vector
 $$X=Y+\sum_{i=1}^{11} c_i X_i$$
 with undetermined  coefficients $c_1,\dots,c _{11}$ to each of the $(n-1)$-dimensional subalgebras (temporarily denoted as  $\hat L$) from the optimal list in \cite{Khab2013}. Then one  determines some of the coefficients by the condition  that $[X,Z_i]\subset\hat L$, where $i=1\dots,{n-1}$ and $Z_1,\dots, Z_{n-1} $ is a basis of $\hat L$. Finally one uses reductions by linear combinations of   $Z_1,\dots, Z_{n-1} $ and those of the inner automorphisms from Table~\ref{table_autom} that stabilize $\hat L$ to  bring $X$ to a ``canonical form'' with as many coefficients as possible being zero. There are two possibilities: either the canonical form of $X$ equals to $Y$ and so $L=\hat L\oplus\{Y\}$  or the canonical form of $X$ can be written as $Z_{n}+Y$, where $Z_n\in L_{11}$ is linearly independent of  $Z_1,\dots, Z_{n-1}$. The first case we call trivial, while in the second case  the subalgebra  with the basis $Z_1,\dots, Z_{n-1}, Z_n+Y$ is included in the optimal list of the $L_{12}$ subalgebras.

\section{Four--dimensional subalgebras and their invariants}\label{sect-inv}
The optimal list of $4$-dimensional subalgebras of $L_{12}$ obtained  in \cite{Sir2014} has 77 items. Most of the items contain a single subalgebra, while some contain parametric families of non-similar subalgebras.  In \cite{Sir2024-partially}, generating sets of  invariants for 29  items  were calculated.  The current paper  treats  20 more items, while   28 more remain to be computed in future works.  Invariants were computed with the help of {\sc Maple} software. In certain cases, computations are simplified by using cylindrical or spherical coordinates. We know give details on coordinate functions and additional variables that are used to write invariants in each of the coordinate systems. Cylindrical and spherical coordinate expressions for the gas dynamic system and for the infinitesimal generators of symmetries listed in Table~\ref{table_sym},  see \cite{SirYulm2020}.

{\bf Cartesian coordinates} are denoted by {\bf (D)} for Descartes and use the following coordinate  functions for  the position and velocity vector 
$$\x=(x, y, z) \text{ and }  \u=(u, v, w)$$
Similar as it has been done in \cite{Khab2019}, in  some cases, one gets  more compact formulas for the invariants, if   $v$ and $w$ are expressed  in terms of $t,y,z$ and  new variables $\bar q$ and $\bar\vartheta $  using change of variables formulas, with $b$ taking various real values:   
\begin{equation}\label{eq1.5.6}   
 v =\frac{ty+bz}{t^2+b^2} + \bar{q} \cos\bar{\vartheta}, \ \ 
 w = \frac{tz-by}{t^2+b^2} + \bar{q} \sin\bar{\vartheta}. 
\end{equation}
If variable change \eqref{eq1.5.6} is applied, it is noted after {\bf (D)}.

{\bf Cylindrical coordinates} are denoted by {\bf (C)}  and use the following coordinate  functions for  the position and velocity vector 
$$\x=(x, r, \thc), \text{ where }r \geqslant 0, \quad 0 \leqslant \thc < 2\pi, \text{ and }  \u=(u, V, W)$$

with the following formulas relating them to their Cartesian counterparts:
\begin{equation}\label{cyl1}y = r \cos \thc, \quad z = r \sin \thc, \quad  v = V \cos \thc - W \sin \thc, \quad w = V \sin \thc + W \cos \thc\end{equation}
It turns out that one gets  more compact formulas for the invariants if   $V$ and $W$ are expressed  in terms of  their own cylindrical coordinates $\qc\geq 0$ and $\varthc\in[0,2\pi)$:
\begin{equation}\label{cyl2}
 V ={\qc} \cos{\varthc}, \ \ 
 W = {\qc} \sin{\varthc}.   
\end{equation}
As a side remark, we note that if we substitute \eqref{cyl2} into the last two equations of \eqref{cyl1}, we get
$$v=\qc\cos(\thc+\varthc) \text{ and  } w=\qc\sin(\thc+\varthc)$$
and observe that $v^2+w^2=V^2+W^2=\qc^2.$

{\bf Spherical coordinates} are denoted by {\bf (S)} and use the following coordinate  functions for  the position and velocity vector 
$$\x=(r_S, \theta_S,\phi), \text{ where }r \geqslant 0, \quad 0 \leqslant \theta_S < \pi, \quad 0 \leqslant \varphi < 2\pi \text{ and }  \u=(U_S, V_S, W_S)$$
with the following formulas relating them to their Cartesian counterparts:
\begin{eqnarray}\label{sph1}& x = r_S \sin \theta_S \cos \phi, \quad y = r_S \sin \theta_S \sin \phi, \quad z = r_S \cos \theta_S,&\\[2mm]
\label{sph2}&u = \left( U_S \sin \theta_S + V_S \cos \theta_S \right) \cos \phi - W_S \sin \phi, \quad
v = \left( U_S \sin \theta_S + V_S \cos \theta_S \right) \sin \phi + W_S \cos \phi, \, &\\[2mm]
&w = U_S \cos \theta_S - V_S \sin \theta_S.& \nonumber
\end{eqnarray}
It turns out that  one gets  more compact formulas for the invariants if $U$,  $V$ and $W$ are expressed  in terms of their own  spherical coordinates $q_S\geq 0$, $\vartheta_S\in [0,\pi)$ and $\varphi \in [0,2\pi) $:
\begin{equation}\label{sph3}
U_S=q_S\,\cos \vartheta_S,  \ \
V_S= q_S\,\sin \vartheta_S  \cos{\varphi}, \ \ 
 W_S = {q_S} \sin{\vartheta_S} \sin{\varphi}.   
\end{equation}
As the a side remark we note that, from \eqref{sph2}-\eqref{sph3} we have 
$$u^2+v^2+w^2=U_S^2+V_S^2+W_S^2=q_S^2$$
and so $q_S$ represents the length of the velocity vector $\u$.
From \eqref{sph1}-\eqref{sph3}, it follows that the scalar product  
$$\x\cdot \u=r_SU_S=r_Sq_S\cos\vartheta_S$$
and therefore $\vartheta_S$ represents  the angle between the velocity and the position vectors.

From the dimensional considerations, we know that a minimal generating set of invariants for each of the four-dimensional subalgebras consists of five invariants.  The density $\rho$ is an invariant for any of the subalgebras and is not listed.  The number of the subalgebra corresponds to the optimal system of subalgebras in \cite{Sir2014}.

\begin{enumerate}[label={}, align=left, leftmargin=0.5cm, itemsep=1em]

\item
\begin{tabular}{@{}p{0.5cm} p{15.2cm}@{}}
\text{4.1 } & \textbf{Basis: } $X_7,\ X_8,\ X_9,\ \Y + X_{11}$ \\
& \textbf{Inv. (S): } $\dfrac{r_S}{t},\ q_S,\ \vartheta_S,\ P - \ln|t|$
\end{tabular}

\item
\begin{tabular}{@{}p{0.5cm} p{15.2cm}@{}}
\text{4.2 } & \textbf{Basis: } $X_7,\ X_8,\ X_9,\ \Y + X_{10}$ \\
& \textbf{Inv. (S): } $r_S,\ q_S,\ \vartheta_S,\ P - t$
\end{tabular}

\item
\begin{tabular}{@{}p{0.5cm} p{15.2cm}@{}}
\text{4.3 } & \textbf{Basis: } $X_1,\ a X_4 + X_7,\ b X_4 + X_{11},\ \Y + X_4$ \\
& \textbf{Inv. (C): } $\dfrac{r}{t},\ \qc,\ \varthc,\ u - P - a \thc - b \ln|t|$
\end{tabular}

\item
\begin{tabular}{@{}p{0.5cm} p{15.2cm}@{}}
\text{4.21 } & \textbf{Basis: } $X_1,\ X_4,\ X_7 + X_{10},\ \Y + X_{10}$ \\
& \textbf{Inv. (C): } $r,\ \qc,\ \varthc,\ P - t + \thc$
\end{tabular}

\item
\begin{tabular}{@{}p{0.5cm} p{15.2cm}@{}}
\text{4.23 } & \textbf{Basis: } $X_2,\ X_3,\ X_7 + X_{10},\ \Y + a X_1 + b X_{10},$  \\
{\centering \text{i.} \par} & $a \neq 0,\,  { a^2 + b^2 = 1}$, \textbf{Inv. (C): } $\dfrac{b}{a} x - t + \thc + \varthc,\ u,\ \qc,\ P - \dfrac{x}{a}$ \\
{\centering \text{ii.} \par} & $a = 0, b = 1$, \qquad\quad \textbf{Inv. (C): } $x,\ u,\ \qc,\ P - t + \thc + \varthc$
\end{tabular}

\item
\begin{tabular}{@{}p{0.5cm} p{15.2cm}@{}}
\text{4.27 } & \textbf{Basis: } $X_1,\ X_4,\ X_{11},\ \Y + X_7$ \\
& \textbf{Inv. (C): } $\dfrac{r}{t},\ \qc,\ \varthc,\ P - \thc$
\end{tabular}

\item
\begin{tabular}{@{}p{0.5cm} p{15.2cm}@{}}
\text{4.34 } & \textbf{Basis: } $X_1,\ X_4,\ X_{10},\ \Y + a X_7 + X_{11}$ \\
{\centering \text{i.} \par} & $a \neq 0$, \textbf{Inv. (C): } $\thc - a \ln|r|,\ \qc,\ \varthc,\ P - \ln|r|$ \\
{\centering \text{ii.} \par} & $a = 0$, \textbf{Inv. (D): } $\dfrac{y}{z},\ v,\ w,\ P - \ln|y|$
\end{tabular}

\item
\begin{tabular}{@{}p{0.5cm} p{15.2cm}@{}}
\text{4.35 } & \textbf{Basis: } $X_1,\ X_4,\ X_{10},\ \Y + X_7$ \\
& \textbf{Inv. (C): } $r,\ \varthc,\ \qc,\ P - \thc$
\end{tabular}

\item
\begin{tabular}{@{}p{0.5cm} p{15.2cm}@{}}
\text{4.38 } & \textbf{Basis: } $X_2,\ X_3,\ \varepsilon X_4 + X_{10},\ \Y + a X_1 + X_7,$\,$ \, {\varepsilon\in\{0,1\} }$ \\
& \textbf{Inv. (C): } $x - \varepsilon \dfrac{t^2}{2} - a (\thc + \varthc),\ u - \varepsilon t,\ \qc,\ P - (\thc + \varthc)$
\end{tabular}

\item
\begin{tabular}{@{}p{0.5cm} p{15.2cm}@{}}
\text{4.42 } & \textbf{Basis: } $a X_1 + X_4,\ b X_3 + X_5,\ b X_2 - X_6,\ \Y + \varepsilon X_1 + X_7,\, a^2+b^2=1, {\varepsilon\in\{0,1\}}$,\,  \\[2mm]
{\centering \par} & $ $ \textbf{Inv. (D),  \eqref{eq1.5.6}:} $t,\ u-\dfrac{x}{t+a} + \dfrac{\varepsilon}{t+a}\bar{\vartheta},\ \bar{q},\ P - \bar{\vartheta}$
\end{tabular}

\item
\begin{tabular}{@{}p{0.5cm} p{15.2cm}@{}}
\text{4.44 } & \textbf{Basis: } $X_4,\ X_5,\ X_6,\ \Y + a X_7 + X_{11}$ \\
{\centering \text{i.} \par} & $a \neq 0$, \textbf{Inv. (D), \eqref{eq1.5.6} with b=0: } $u - \dfrac{x}{t},\ \ln|t| - \dfrac{\bar{\vartheta}}{a},\ \bar{q},\ P - \ln|t|$ \\
{\centering \text{ii.} \par} & $a = 0$, \textbf{Inv. (D): } $u - \dfrac{x}{t},\ v - \dfrac{y}{t},\ w - \dfrac{z}{t},\ P - \ln|t|$
\end{tabular}

\item
\begin{tabular}{@{}p{0.5cm} p{15.2cm}@{}}
\text{4.45 } & \textbf{Basis: } $X_4,\ X_5,\ X_6,\ \Y + \varepsilon X_1 + X_7,\,\, { \varepsilon\in\{0,1\}}\,$ \\
&  \textbf{Inv. (D), \eqref{eq1.5.6} with b=0:} $t,\ u - \dfrac{x}{t} + \dfrac{\varepsilon}{t} \bar{\vartheta},\ \bar{q},\ P - \bar{\vartheta}$
\end{tabular}

\item
\begin{tabular}{@{}p{0.5cm} p{15.2cm}@{}}
\text{4.54 } & \textbf{Basis: } $X_1,\ X_3 + X_5,\ X_2 - X_6,\ \Y + a X_4 + X_7$ \\
& \textbf{Inv. (D), \eqref{eq1.5.6} with b=1:} $t,\ u - a \bar{\vartheta},\ \bar{q},\ P - \bar{\vartheta}$
\end{tabular}

\item
\begin{tabular}{@{}p{0.5cm} p{15.2cm}@{}}
\text{4.56 } & \textbf{Basis: } $X_1,\ X_5,\ X_6,\ \Y + a X_4 + b X_7 + X_{11}$ \\
{\centering \text{i.} \par} & $b \neq 0$, \textbf{Inv. (D), \eqref{eq1.5.6} with b=0: } $u - a\ln|t|,\ \bar{\vartheta} - b \ln|t|,\ \bar{q},\ P - \ln|t|$ \\
{\centering \text{ii.} \par} & $b = 0$, \textbf{Inv. (D): } $u - a \ln|t|,\ v - \dfrac{y}{t},\ w - \dfrac{z}{t},\ P - \ln|t|$
\end{tabular}

\item
\begin{tabular}{@{}p{0.5cm} p{15.2cm}@{}}
\text{4.57 } & \textbf{Basis: } $X_1,\ X_5,\ X_6,\ \Y + b X_4 + X_7$ \\
&  \textbf{Inv. (D), \eqref{eq1.5.6} with b=0: } $t,\ u - b \bar{\vartheta},\ \bar{q},\ P - \bar{\vartheta}$
\end{tabular}

\item
\begin{tabular}{@{}p{0.5cm} p{15.2cm}@{}}
\text{4.64 } & \textbf{Basis: } $X_2,\ X_3,\ X_4,\ \Y + a X_7 + X_{11}$ \\
{\centering \text{i.} \par} & $a \neq 0$, \textbf{Inv. (C): } $u - \dfrac{x}{t},\ \thc + \varthc - a \ln|t|,\ \qc,\ P - \ln|t|$ \\
{\centering \text{ii.} \par} & $a = 0$, \textbf{Inv. (D): } $u - \dfrac{x}{t},\ v,\ w,\ P - \ln|t|$
\end{tabular}

\item
\begin{tabular}{@{}p{0.5cm} p{15.2cm}@{}}
\text{4.65 } & \textbf{Basis: } $X_2,\ X_3,\ X_4,\ \Y + \varepsilon X_1 + X_7, \, { \varepsilon\in\{0,1\}}\,$ \\
& \textbf{Inv. (C): } $t,\ \qc,\ u - \dfrac{x}{t} + \dfrac{\varepsilon}{t} (\thc + \varthc),\ P - (\thc + \varthc)$
\end{tabular}

\item
\begin{tabular}{@{}p{0.5cm} p{15.2cm}@{}}
\text{4.71} & \textbf{Basis: } $X_1,\ X_2,\ X_4,\ \Y + a X_5 + b X_6 + c X_{10} + d X_{11} $ \\[2mm]
{\centering \text{i.} \par} & $d\neq 0$, $c^2+d^2=1$, \textbf{Inv. (D): }$\dfrac{z}{dt + c} - \dfrac{bc}{d^2 (dt + c)} - \dfrac{b \ln|dt + c|}{d^2}$,\\  & $v - \dfrac{a}{d} \ln|dt + c|, \ \ w - \dfrac{b}{d} \ln|dt + c|, \ \ P - \dfrac{\ln|dt + c|}{d}$ \\[2mm]
{\centering \text{ii.} \par} & $d=0, c=1$, \textbf{Inv. (D): }$z-b\dfrac{t^2}{2}$,\ \ $v-at$,\ \ $w-bt$,\ \ $P-t$ 
\end{tabular}

\item
\begin{tabular}{@{}p{0.5cm} p{15.2cm}@{}}
\text{4.74 } & \textbf{Basis: } $X_1,\ X_2,\ X_3,\ \Y + a X_4 + X_7 + \varepsilon X_{10} + c X_{11}$  \\
{\centering \text{i.} \par} & $c \neq 0, \varepsilon = 0$, \textbf{Inv. (C): } $\thc + \varthc - \dfrac{\ln|t|}{c},\ u - \dfrac{a}{c} \ln|t|,\ \qc,\ P - \dfrac{\ln|t|}{c}$ \\
{\centering \text{ii.} \par} & $c = 0, \varepsilon = 0$, \textbf{Inv. (C): } $t,\ U - a (\varthc + \thc),\ \qc,\ P - (\thc + \varthc)$ \\
{\centering \text{iii.} \par} & $c = 0, \varepsilon = 1$, \textbf{Inv. (C): } $\thc + \varthc - t,\ u - a t,\ \qc,\ P - t$
\end{tabular}

\item
\begin{tabular}{@{}p{0.5cm} p{15.2cm}@{}}
\text{4.77 } & \textbf{Basis: } $X_1,\ X_2,\ X_3,\ \Y + X_4$ \\
& \textbf{Inv. (D): } $t,\ v,\ w,\ P - u$
\end{tabular}

\end{enumerate}

\setlength\extrarowheight{8pt}

\section{Isomorphism classes}\label{sect-isom}
This section is devoted to the classification, up to a Lie algebra isomorphism,  of the four-dimensional non-similar subalgebras listed in the previous section. The isomorphism classes of four-dimensional real Lie algebras are listed  in~\cite{PatWint1976} and~\cite{PatWint1977}. In combination with the classification of sub-algebra inclusions in~\cite{PatWint1977}, the results of this section provide a seed for future   studies of the submodel hierarchy and the corresponding nesting of invariant and partially invariant solutions.

In~\cite{PatWint1976} and~\cite{PatWint1977},
notation $A_{r,j}$ is used to denote an indecomposable Lie algebra of  dimension $r$, whose position in their classification  list  is $j$. The superscript(s), if any, denote(s) continuous parameter(s) on which the algebra depends.   There is a single isomorphism class $A_1$ of dimension one. There are two isomorphism classes of dimension two, abelian $2A_1=A_1\oplus A_1$ and solvable   $A_2$. Isomorphism classes  of three-dimensional indecomposable Lie algebras are presented by nine types, some with parameters,   denoted  $A_{3,i}$ (with $i = 1, ... ,9$). The are two isomorphism classes of decomposable three-dimensional Lie algebras: abelian $3A_1$ and non-abelian $A_2\oplus A_1$. Isomorphism classes  of four-dimensional indecomposable Lie algebras are presented  by twelve types $A_{4,i}$ (with $i = 1, ... ,12$), some with parameters. The are twelve isomorphism classes of decomposable  four-dimensional Lie algebras: $4A_1$, $2A_2$, $A_2\oplus 2A_1$, and  $A_1\oplus A_{3,i}$, $(i=1, ... ,9)$.

The following conventions are used in the list below. The first column indicates the subalgebra number \textbf{N} from the list in the previous section. The second column indicates  the isomorphism class presented in~\cite{PatWint1976} and~\cite{PatWint1977}, with the convention described in the previous paragraph. The third column gives a change of basis formulas, where $E_1$, $E_2$, $E_3$, $E_4$ are the basis infinitesimal generators listed in the previous section.  Lastly, the fourth column indicates nonzero commutator relations for the new basis, which match the commutator relations from the table in~\cite{PatWint1977}. 

\begin{enumerate}[label={}, align=left, leftmargin=0cm, itemsep=1em]

\item \textbf{\parbox[t]{1.2cm}{N} \parbox[t]{3.3cm}{\quad Algebra} \parbox[t]{4cm}{Change of the basis} \parbox[t]{4.92cm}{Nonzero commutators}}

\item \parbox[t]{1.2cm}{4.1} \parbox[t]{3.3cm}{$A_{3,9} \oplus A_1$} \parbox[t]{4cm}{$e_1=E_1,\ \ e_2=-E_2,\\[2mm] e_3=E_3,\ \ e_4=E_4$} \parbox[t]{4.92cm}{$\left[e_1, e_2\right] = e_3$,\\[2mm] $\left[e_2, e_3\right] = e_1$,\\[2mm] $\left[e_3, e_1\right] = e_2$}

\item \parbox[t]{1.2cm}{4.2} \parbox[t]{3.3cm}{$A_{3,9} \oplus A_1$} \parbox[t]{4cm}{$e_1=E_1,\ \ e_2=-E_2,\\[2mm] e_3=E_3,\ \ e_4=E_4$} \parbox[t]{4.92cm}{$\left[e_1, e_2\right] = e_3$,\\[2mm] $\left[e_2, e_3\right] = e_1$,\\[2mm] $\left[e_3, e_1\right] = e_2$}

\item \parbox[t]{1.2cm}{4.3} \parbox[t]{3.3cm}{$A_2 \oplus 2A_1$} \parbox[t]{4cm}{$e_1=-E_3,\ \ e_2=E_1,\\[2mm] e_3=E_2,\ \ e_4=E_4$} \parbox[t]{4.92cm}{$\left[e_1, e_2\right] = e_2$}

\item \parbox[t]{1.2cm}{4.21} \parbox[t]{3.3cm}{$A_{3,1} \oplus A_1$} \parbox[t]{4cm}{$e_1=-E_1,\ \ e_2=E_2,\\[2mm] e_3=E_3,\ \ e_4=E_4$} \parbox[t]{4.92cm}{$\left[e_2, e_3\right] = e_1$}

\item \parbox[t]{1.2cm}{4.23 i} \parbox[t]{3.3cm}{$A_{3,6} \oplus A_1$} \parbox[t]{4cm}{$e_1=-E_1,\ \ e_2=E_2,\\[2mm] e_3=E_3,\ \ e_4=E_4$} \parbox[t]{4.92cm}{$\left[e_1, e_3\right] = -e_2$, $\left[e_2, e_3\right] = e_1$}

\item \parbox[t]{1.2cm}{4.23 ii} \parbox[t]{3.3cm}{$A_{3,6} \oplus A_1$} \parbox[t]{4cm}{$e_1=-E_1,\ \ e_2=E_2,\\[2mm] e_3=E_3,\ \ e_4=E_4$} \parbox[t]{4.92cm}{$\left[e_1, e_3\right] = -e_2$, $\left[e_2, e_3\right] = e_1$}

\item \parbox[t]{1.2cm}{4.27} \parbox[t]{3.3cm}{$A_2 \oplus 2A_1$} \parbox[t]{4cm}{$e_1=-E_3,\ \ e_2=E_1,\\[2mm] e_3=E_2,\ \ e_4=E_4$} \parbox[t]{4.92cm}{$\left[e_1, e_2\right] = e_2$}

\item \parbox[t]{1.2cm}{4.34 i} \parbox[t]{3.3cm}{$A_{4,9}^0$} \parbox[t]{4cm}{$e_1=E_1,\ \ e_2=E_3,\\[2mm] e_3=E_2,\ \ e_4=E_4$} \parbox[t]{4.92cm}{$\left[e_2, e_3\right] = e_1$, $\left[e_1, e_4\right] = e_1$,\\[2mm] $\left[e_2, e_4\right] = e_2$}

\item \parbox[t]{1.2cm}{4.34 ii} \parbox[t]{3.3cm}{$A_{4,9}^0$} \parbox[t]{4cm}{$e_1=E_1,\ \ e_2=E_3,\\[2mm] e_3=E_2,\ \ e_4=E_4$} \parbox[t]{4.92cm}{$\left[e_2, e_3\right] = e_1$,\ \  $\left[e_1, e_4\right] = e_1$,\\[2mm] $\left[e_2, e_4\right] = e_2$}

\item \parbox[t]{1.2cm}{4.35} \parbox[t]{3.3cm}{$A_{3,1} \oplus A_1$} \parbox[t]{4cm}{$e_1=-E_1,\ \ e_2=E_2,\\[2mm] e_3=E_3,\ \ e_4=E_4$} \parbox[t]{4.92cm}{$\left[e_2, e_3\right] = e_1$}

\item \parbox[t]{1.2cm}{4.38} \parbox[t]{3.3cm}{$A_{3,6} \oplus A_1$} \parbox[t]{4cm}{$e_1=E_2,\ \ e_2=E_1,\\[2mm] e_3=E_4,\ \ e_4=E_3$} \parbox[t]{4.92cm}{$\left[e_1, e_3\right] = -e_2$, $\left[e_2, e_3\right] = e_1$}

\item \parbox[t]{1.2cm}{4.42} \parbox[t]{3.3cm}{$A_{3,6} \oplus A_1$} \parbox[t]{4cm}{$e_1=E_2,\ \ e_2=E_3,\\[2mm] e_3=E_4,\ \ e_4=E_1$} \parbox[t]{4.92cm}{$\left[e_1, e_3\right] = -e_2$, $\left[e_2, e_3\right] = e_1$}

\item \parbox[t]{1.2cm}{4.44 i} \parbox[t]{3.3cm}{$A_{3,6} \oplus A_1$} \parbox[t]{4cm}{$e_1=E_3,\ \ e_2=E_2,\\[2mm] e_3=E_4/a,\ \ e_4=E_1$} \parbox[t]{4.92cm}{$\left[e_1, e_3\right] = -e_2$, $\left[e_2, e_3\right] = e_1$}

\item \parbox[t]{1.2cm}{4.44 ii} \parbox[t]{3.3cm}{$4A_1$} \parbox[t]{4cm}{$e_1=E_1,\ \ e_2=E_2,\\[2mm] e_3=E_3,\ \ e_4=E_4$} \parbox[t]{4.92cm}{ }

\item \parbox[t]{1.2cm}{4.45} \parbox[t]{3.3cm}{$A_{3,6} \oplus A_1$} \parbox[t]{4cm}{$e_1=E_3,\ \ e_2=E_2,\\[2mm] e_3=E_4,\ \ e_4=E_1$} \parbox[t]{4.92cm}{$\left[e_1, e_3\right] = -e_2$, $\left[e_2, e_3\right] = e_1$}

\item \parbox[t]{1.2cm}{4.54} \parbox[t]{3.3cm}{$A_{3,6} \oplus A_1$} \parbox[t]{4cm}{$e_1=E_2,\ \ e_2=E_3,\\[2mm] e_3=E_4,\ \ e_4=E_1$} \parbox[t]{4.92cm}{$\left[e_1, e_3\right] = -e_2$, $\left[e_2, e_3\right] = e_1$}

\item \parbox[t]{1.2cm}{4.56 i} \parbox[t]{3.3cm}{$A_{4,6}^{\frac{1}{|b|},0}, b\neq 0$} \parbox[t]{4cm}{$e_1=E_1,\ \ e_2=-\frac {b}{|b|}E_2,\\[2mm] e_3=E_3,\ \ e_4=\frac 1 {|b|}E_4$} \parbox[t]{4.92cm}{$\left[e_1, e_4\right] = \frac 1 {|b|} e_1$,\\[2mm] $\left[e_2, e_4\right] = -e_3$,\\[2mm] $\left[e_3, e_4\right] = e_2$}

\item \parbox[t]{1.2cm}{4.56 ii} \parbox[t]{3.3cm}{$A_2 \oplus 2A_1$} \parbox[t]{4cm}{$e_1=-E_4,\ \ e_2=E_1,\\[2mm] e_3=E_3,\ \ e_4=E_2$} \parbox[t]{4.92cm}{$\left[e_1, e_2\right] = e_2$}

\item \parbox[t]{1.2cm}{4.57} \parbox[t]{3.3cm}{$A_{3,6} \oplus A_1$} \parbox[t]{4cm}{$e_1=E_3,\ \ e_2=E_2,\\[2mm] e_3=E_4,\ \ e_4=E_1$} \parbox[t]{4.92cm}{$\left[e_1, e_3\right] = -e_2$, $\left[e_2, e_3\right] = e_1$}

\item \parbox[t]{1.2cm}{4.64 i} \parbox[t]{3.3cm}{$A_{3,7}^{\frac{1}{|a|}} \oplus A_1$, $a \neq 0$} \parbox[t]{4cm}{$e_1=E_2,\ \ e_2=\frac{a}{|a|}E_1,\\[2mm] e_3=\frac {1}{|a|}E_4,\ \ e_4=E_3$} \parbox[t]{4.92cm}{$\left[e_1, e_3\right] =\frac 1{|a|} e_1 - e_2$,\\[2mm] $\left[e_2, e_3\right] = e_1 + \frac 1 {|a|}e_2$}

\item \parbox[t]{1.2cm}{4.64 ii} \parbox[t]{3.3cm}{$A_{3,3} \oplus A_1$} \parbox[t]{4cm}{$e_1=E_1,\ \ e_2=E_2,\\[2mm] e_3=E_4,\ \ e_4=E_3$} \parbox[t]{4.92cm}{$\left[e_1, e_3\right] = e_1$, $\left[e_2, e_3\right] = e_2$}

\item \parbox[t]{1.2cm}{4.65} \parbox[t]{3.3cm}{$A_{3,6} \oplus A_1$} \parbox[t]{4cm}{$e_1=E_2,\ \ e_2=E_1,\\[2mm] e_3=E_4,\ \ e_4=E_3$} \parbox[t]{4.92cm}{$\left[e_1, e_3\right] = -e_2$, $\left[e_2, e_3\right] = e_1$}

\item \parbox[t]{1.2cm}{4.71 i} \parbox[t]{3.3cm}{$A_{3,3} \oplus A_1$} \parbox[t]{4cm}{$e_1=E_1,\ \ e_2=E_2,\\[2mm] e_3=\frac{1}{d}E_4,\ \ e_4=\frac{c}{d}E_1+E_3$} \parbox[t]{4.92cm}{$\left[e_1, e_3\right] = e_1$, $\left[e_2, e_3\right] = e_2$}

\item \parbox[t]{1.2cm}{4.71 ii} \parbox[t]{3.3cm}{$A_{3,1} \oplus A_1$} \parbox[t]{4cm}{$e_1=E_1,\ \ e_2=E_4,\\[2mm] e_3=E_3,\ \ e_4=E_2$} \parbox[t]{4.92cm}{$\left[e_2, e_3\right] = e_1$}

\item \parbox[t]{1.2cm}{4.74 i} \parbox[t]{3.3cm}{$A_{4,6}^{|c|,|c|}$, $c \neq 0$} \parbox[t]{4cm}{$e_1=E_1,\ \ e_2=E_3,\\[2mm] e_3=\frac{c}{|c|}E_2,\ \ e_4=\frac{c}{|c|}E_4$} \parbox[t]{4.92cm}{$\left[e_1, e_4\right] = |c| e_1$,\\[2mm] $\left[e_2, e_4\right] = |c|e_2 - e_3$,\\[2mm] $\left[e_3, e_4\right] = e_2 + |c|e_3$}

\item \parbox[t]{1.2cm}{4.74 ii} \parbox[t]{3.3cm}{$A_{3,6} \oplus A_1$} \parbox[t]{4cm}{$e_1=E_3,\ \ e_2=E_2,\\[2mm] e_3=E_4,\ \ e_4=E_1$} \parbox[t]{4.92cm}{$\left[e_1, e_3\right] = -e_2$, $\left[e_2, e_3\right] = e_1$}

\item \parbox[t]{1.2cm}{4.74 iii} \parbox[t]{3.3cm}{$A_{3,6} \oplus A_1$} \parbox[t]{4cm}{$e_1=E_3,\ \ e_2=E_2,\\[2mm] e_3=E_4,\ \ e_4=E_1$} \parbox[t]{4.92cm}{$\left[e_1, e_3\right] = -e_2$, $\left[e_2, e_3\right] = e_1$}

\item \parbox[t]{1.2cm}{4.77} \parbox[t]{3.3cm}{$4A_1$} \parbox[t]{4cm}{$e_1=E_1,\ \ e_2=E_2,\\[2mm] e_3=E_3,\ \ e_4=E_4$} \parbox[t]{4.92cm}{}

\end{enumerate}

\section{A submodel and exact solutions}\label{sect-sol}
From the generating sets of invariants listed in Section~\ref{sect-inv}, one can observe that  the corresponding submodels would be either of rank 1 and defect 1 or of rank 0 and defect 0.  Here we construct a submodel   for the  subalgebra with number 77  of rank 1 and defect 1. The basis generators in the Cartesian coordinate system have the form
\begin{equation}
    X_1=\partial_x,\ \ X_2=\partial_y,\ \ X_3=\partial_z,\ \ \Y+X_4=\partial_P+t\partial_x+\partial_u.\label{eqn1}
\end{equation}
The invariants of \eqref{eqn1} are
\begin{equation}
t,\ \ v,\ \ w, \ \ P_1=P-u, \ \ \rho.\label{eqn2}
\end{equation}
Using invariants \eqref{eqn2}, we introduce a new set of invariant unknown  functions of invariant independent variable $t$
\begin{equation}
v=v(t),\ \ w=w(t),\ \ \rho=\rho(t),\ \  P_1=P_1(t)\label{eqn3}
\end{equation}
which along with  a non-invariant ``defective''  unknown function
$$u=u(t,x,y,z),$$
of all "independent variables.
We substitute in~\eqref{eq-gd}
 to obtain a reduced system of PDEs, called a \emph{partially invariant submodel of rank 1 and defect 1:}
\begin{equation}\label{eqn4}
\begin{array}{c}
u_t+uu_x+vu_y+wu_z+\rho^{-1}u_x=0,\\
v_t+\rho^{-1}u_y=0,\\
w_t+\rho^{-1}u_z=0,\\
\rho_t+\rho u_x=0,\\
P_{1t}+u_t+uu_x+vu_y+wu_z+\rho f_{\rho}u_x=0,
\end{array}
\end{equation}
where $f$ is an arbitrary function of one variable from  the state equation \eqref{state-sp}.

\subsection{Solution for isochoric motion of media}
If the density is constant $\rho=\rho_0>0$ the volume remains the same and such solutions are called isochoric.  In this case, the exact solution of system \eqref{eqn4} has the form
\begin{equation}\label{eqn_5}
\begin{array}{c}
u=k_0y+m_0z+\dfrac{k_0^2+m_0^2}{2\rho_0}t^2-(k_0v_0+m_0w_0)t+n_0,\\
v=-\dfrac{k_0}{\rho_0}t+v_0,\\
w=-\dfrac{m_0}{\rho_0}t+w_0,\\
\rho=\rho_0,\\
P_1=P_0,
\end{array}
\end{equation}
where $k_0$, $m_0$, $v_0$, $n_0$, $P_0$ are arbitrary constants. The number of arbitrary constants in the solution can be reduced by Galilean translations (see Table~\ref{table_sym}) and pressure translations (see Table~\ref{table_sym_ext}):
\begin{equation}\label{eqn5}
\begin{array}{c}
u=k_0y+m_0z+\dfrac{k_0^2+m_0^2}{2\rho_0}t^2,\\
v=-\dfrac{k_0}{\rho_0}t,\\
w=-\dfrac{m_0}{\rho_0}t,\\
\rho=\rho_0,\\
P_1=0,
\end{array}
\end{equation}
Since $P_1=P-u$, we conclude that in this case
\begin{equation}\label{Pisoch}P=u=k_0y+m_0z+\dfrac{k_0^2+m_0^2}{2\rho_0}t^2,
\end{equation}
while from the state equation  \eqref{state-sp}, we conclude that
\begin{equation}
    \label{Sisoch} S=P-f(\rho)=k_0y+m_0z+\dfrac{k_0^2+m_0^2}{2\rho_0}t^2-f(\rho_0)
\end{equation}
The first four equations in \eqref{eqn5}, along with \eqref{Pisoch} and \eqref{Sisoch} is an exact solution of  the gas-dynamic system \eqref{eq-gd} with state equation \eqref{state-sp}.

From \eqref{eqn5} we see that  motion of particles  has a vortex 
\begin{equation}\label{eqn21_1}
	\begin{array}{c}
		\text{rot}\,\u=(w_y-v_z,u_z-w_x,v_x-u_y)=\\
		=\left(0,m_0,-k_0\right).
	\end{array}
\end{equation}

The position of particles is described by the equation:
\begin{equation}\label{eqn6}
\dfrac{d\x}{dt}=\u(t,\x).
\end{equation}
Solving \eqref{eqn6} with the right-hand side prescribed \eqref{eqn5} we obtain equations of the particle flow
\begin{equation}\label{eqn20} 
	\begin{array}{c}
		x(t)=(k_0y_0+m_0z_0)t+x_0,\\
		y(t)=-\dfrac{k_0}{2\rho_0}t^2+y_0,\\
		z(t)=-\dfrac{m_0}{2\rho_0}t^2+z_0,
	\end{array}
\end{equation}
  where $x_0, y_0, z_0$ are the coordinates of the particle at $t=0$.

\begin{figure}[!h]
	\centering
	\includegraphics[scale=0.4]{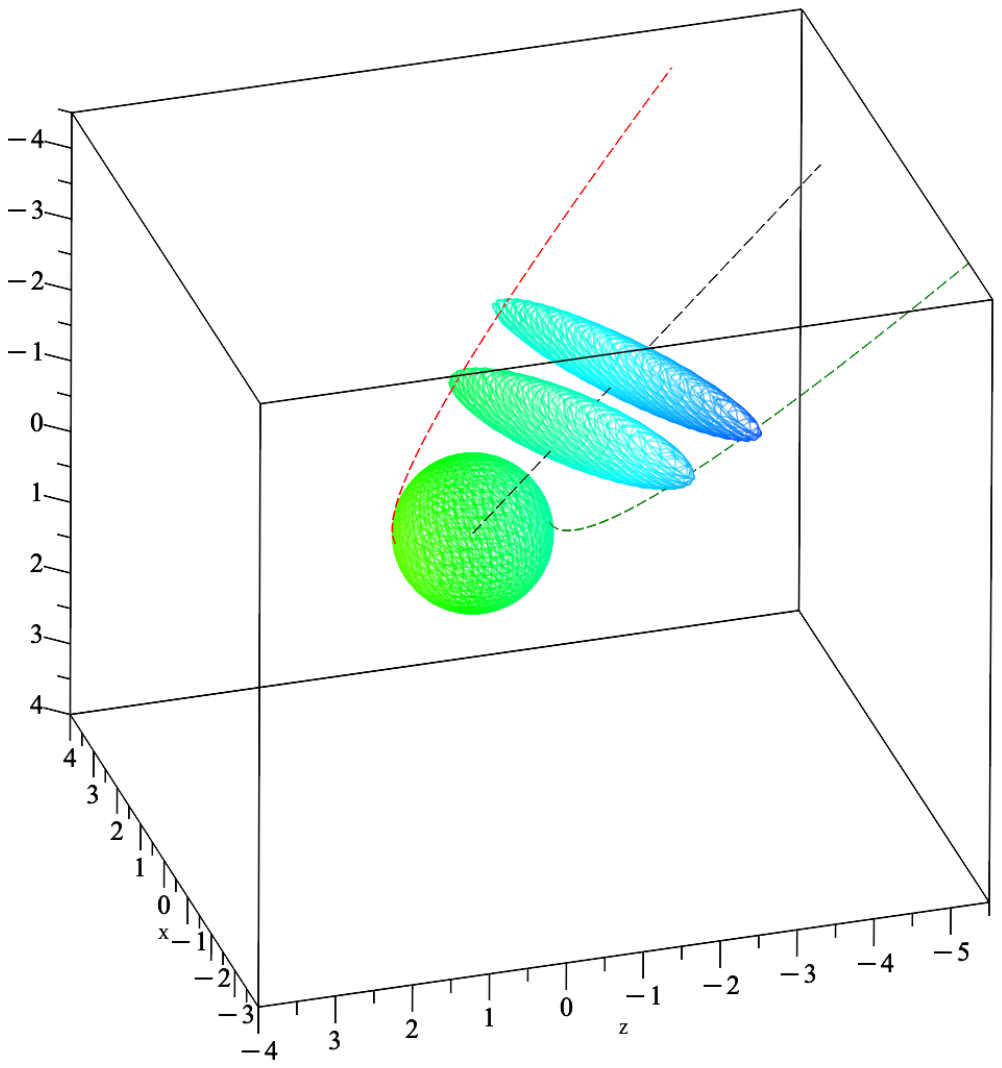}
	\caption{The motion of the particle volume~\eqref{eqn8},~\eqref{eqn20} with $\rho_0=1$, $k_0=1$, $m_0=1$. The volume is a sphere at $t=0$ and becomes ellipsoidal at $t=1.6$, $t=2$. The trajectories of particles~\eqref{eqn20} with initial coordinates $(x_0,y_0,z_0)=(0,0,0)$ (black dash curve); $(0,0,1)$ (red dash curve); $(0,0,-1)$ (green dash curve) are shown for $t=0$ to $t=3$.}\label{fig1}
\end{figure}

One can think ~\eqref{eqn20} as time-parametrized map from $\mathbb R^3$ to $\mathbb R^3$  from the initial position to the position at time $t$ and observe that its Jacobian determinant is constant confirming that the map is volume preserving as expected: 
\begin{equation}
	\begin{array}{c}
		J=\begin{vmatrix}
			1& k_0t & m_0t\\
			0& 1& 0\\
			0& 0& 1
		\end{vmatrix}
	\end{array}=1.
\end{equation}

Let particles \eqref{eqn20} be on the sphere with radius being equal to 1 at the moment of time  $t=0$ 
\begin{equation}\label{eqn8}
	\begin{array}{c}
		x_0^2+y_0^2+z_0^2=1.
	\end{array}
\end{equation}
Over time, the particles will form an ellipsoid, with their volumes coincident (Fig.~\ref{fig1}).

Equations \eqref{eqn20} can be viewed as a transition from Eulerian coordinates (with time and \emph{current} particle position being independent variables) to Lagrangian coordinates (with time and \emph{initial} particle position being independent variables). 
In the Lagrangian coordinates  the components of the velocity vector  obtained by the taking time derivative of \eqref{eqn20}:
\begin{equation}
u=k_0y_0+m_0z_0,\ \ v=-\dfrac{k_0}{\rho_0}t,\ \ w=-\frac{m_0}{\rho_0}t. \label{eqn0}
\end{equation}
On the other hand,  after substituting \eqref{eqn20} into \eqref{Pisoch}  and \eqref{Sisoch}, we can write the other three gas dynamic functions in terms of the Lagrangian coordinates
\begin{equation}
 \rho=\rho_0,\ \ P=k_0y_0+m_0z_0, \ \ S=k_0y_0+m_0z_0-f(\rho_0). \label{rhoP}   
\end{equation}

Differentiating \eqref{eqn0} we obtain  the acceleration vector: $$\a=\left(0,-\dfrac{k_0}{\rho_0},-\dfrac{m_0}{\rho_0}\right).$$
The particle moves along the $x$-axis with constant velocity. Along the $y$ and $z$-axes, the motion of the particle has constant acceleration.

Projections of the world lines \eqref{eqn20} to $\mathbb{R}^3$ parametrized by the spatial coordinates $(x,y,z)$ are particle trajectories. 
When $k_0=0$ and $m_0=0$ then all particles are stationary.

If $k_0y_0+m_0z_0=0$,  equations~\eqref{eqn20} show then the particle's trajectory  is  a straight line
{
\begin{equation*}
{m_0}(y-y_0)-k_0(z-z_0)=0.
\end{equation*}
}
in the plane $x=x_0$.
In this case \eqref{rhoP}  shows that the pressure $P=0$.

If $k_0=0$, $m_0\neq 0$ equations~\eqref{eqn20} show that a trajectory of a particle  is a parabola
\begin{equation*}
z=z_0-\dfrac{(x-x_0)^2}{2\rho_0m_0z_0^2}.
\end{equation*}
 lying in the plane  $y=y_0$.

If $k_0y_0+m_0z_0\neq 0$ and  $k_0\neq 0$ the projections of a trajectory onto the $(y,z)$-plane is a straight line with a slope of $\dfrac{m_0}{k_0}$, 
\begin{equation*}
z=\dfrac{m_0}{k_0}(y-y_0)+z_0.
\end{equation*}
The projection onto the $(x,y)$-plane is a parabola:
\begin{equation*}
y=y_0-\dfrac{k_0(x-x_0)^2}{2\rho_0(k_0y_0+m_0z_0)^2}.
\end{equation*}
The vertex of the parabola is at the point $(x_0,y_0)$.

The projection onto the $(x,z)$-plane is also a parabola:
\begin{equation*}
z=z_0-\dfrac{m_0(x-x_0)^2}{2\rho_0(k_0y_0+m_0z_0)^2},
\end{equation*}
whose vertex is at the point $(x_0,z_0)$.

\subsection{Solution for non-isochoric motion of media}
If the density $\rho$ is not a constant, exact solution of system~\eqref{eqn4} has the form
\begin{equation}\label{eqn5_1_}
\begin{array}{c}
u=\dfrac{x}{t}+k_0\dfrac{y}{t}+m_0\dfrac{z}{t}+\dfrac{n_0}{t}+\dfrac{k_0^2+m_0^2-1}{2\rho_0}t-k_0v_0-m_0w_0,\\
v=-\dfrac{k_0}{\rho_0}t+v_0,\\
w=-\dfrac{m_0}{\rho_0}t+w_0,\\
\rho=\dfrac{\rho_0}{t},\\
P_1=f\left(\dfrac{\rho_0}{t}\right)+\dfrac{t}{\rho_0}+P_0.
\end{array}
\end{equation}
where $k_0$, $m_0$, $v_0$, $n_0$, $w_0$, $P_0$   and $\rho_0>0$ are arbitrary constants, and the number of arbitrary constants in the solution has been reduced by space translations, Galilean translations (see Table~\ref{table_sym}) and pressure translations (see Table~\ref{table_sym_ext})
\begin{equation}\label{eqn5_1}
\begin{array}{c}
u=\dfrac{x}{t}+k_0\dfrac{y}{t}+m_0\dfrac{z}{t}+\dfrac{k_0^2+m_0^2-1}{2\rho_0}t,\\
v=-\dfrac{k_0}{\rho_0}t,\\
w=-\dfrac{m_0}{\rho_0}t,\\
\rho=\dfrac{\rho_0}{t},\\
P_1=f\left(\dfrac{\rho_0}{t}\right)+\dfrac{t}{\rho_0}.
\end{array}
\end{equation}
Since $P_1=P-u$, we conclude that in this case
\begin{equation}
    \label{Pnisoch}P=\dfrac{x}{t}+k_0\dfrac{y}{t}+m_0\dfrac{z}{t}+\dfrac{k_0^2+m_0^2-1}{2\rho_0}t+f\left(\dfrac{\rho_0}{t}\right)+\dfrac{t}{\rho_0},
\end{equation}
while from the state equation  \eqref{state-sp}, we conclude that
\begin{equation}
    \label{Snisoch} S=P-f(\rho)=\dfrac{x}{t}+k_0\dfrac{y}{t}+m_0\dfrac{z}{t}+\dfrac{k_0^2+m_0^2-1}{2\rho_0}t+\dfrac{t}{\rho_0}.
    \end{equation}
    The first four equations in \eqref{eqn5_1}, along with \eqref{Pnisoch} and \eqref{Snisoch} is an exact solution of  the gas-dynamic system \eqref{eq-gd} with state equation \eqref{state-sp}.

From \eqref{eqn5_1} we see that  motion of particles  has a vortex 
\begin{equation*}
		\text{rot}\,\u=(w_y-v_z,u_z-w_x,v_x-u_y)
		=\left(0,\dfrac{m_0}{t},-\dfrac{k_0}{t}\right).
\end{equation*}
    
The world lines of the particles in $\mathbb{R}^4(t,\x)$, obtained  from \eqref{eqn6} and \eqref{eqn5_1}, are 
\begin{equation}\label{eqn7} 
	\begin{array}{c}
		x=-(k_0y_0+m_0z_0)-\dfrac{t^2}{2\rho_0}+u_0t,\\
		y=-\dfrac{k_0}{2\rho_0}t^2+y_0,\\
		z=-\dfrac{m_0}{2\rho_0}t^2+z_0
	\end{array}
\end{equation}
and can be viewed  as a transition from  the Eulerian  to Lagrangian coordinates. In the Lagrangian coordinates  the components of the velocity vector  obtained by the taking time derivative of \eqref{eqn7}:
\begin{equation}\label{eqn7-1} 
	\begin{array}{c}
		u=-\dfrac{t}{\rho_0}+u_0,\\
		v=-\dfrac{k_0}{\rho_0}t,\\
		w=-\dfrac{m_0}{\rho_0}t.
	\end{array}
\end{equation}
while the other three gas-dynamic variables  become
$$ \rho=\dfrac{\rho_0}{t},\ \ P=u_0+f\left(\dfrac{\rho_0}{t}\right),\ \ S=u_0.$$
Taking the time-derivative of \eqref{eqn7-1}, we see that particles are moving with  the constant acceleration:
$$\a=\left(-\dfrac{1}{\rho_0},-\dfrac{k_0}{\rho_0},-\dfrac{m_0}{\rho_0}\right).$$
\begin{figure}[!h]
	\centering
	\includegraphics[scale=0.4]{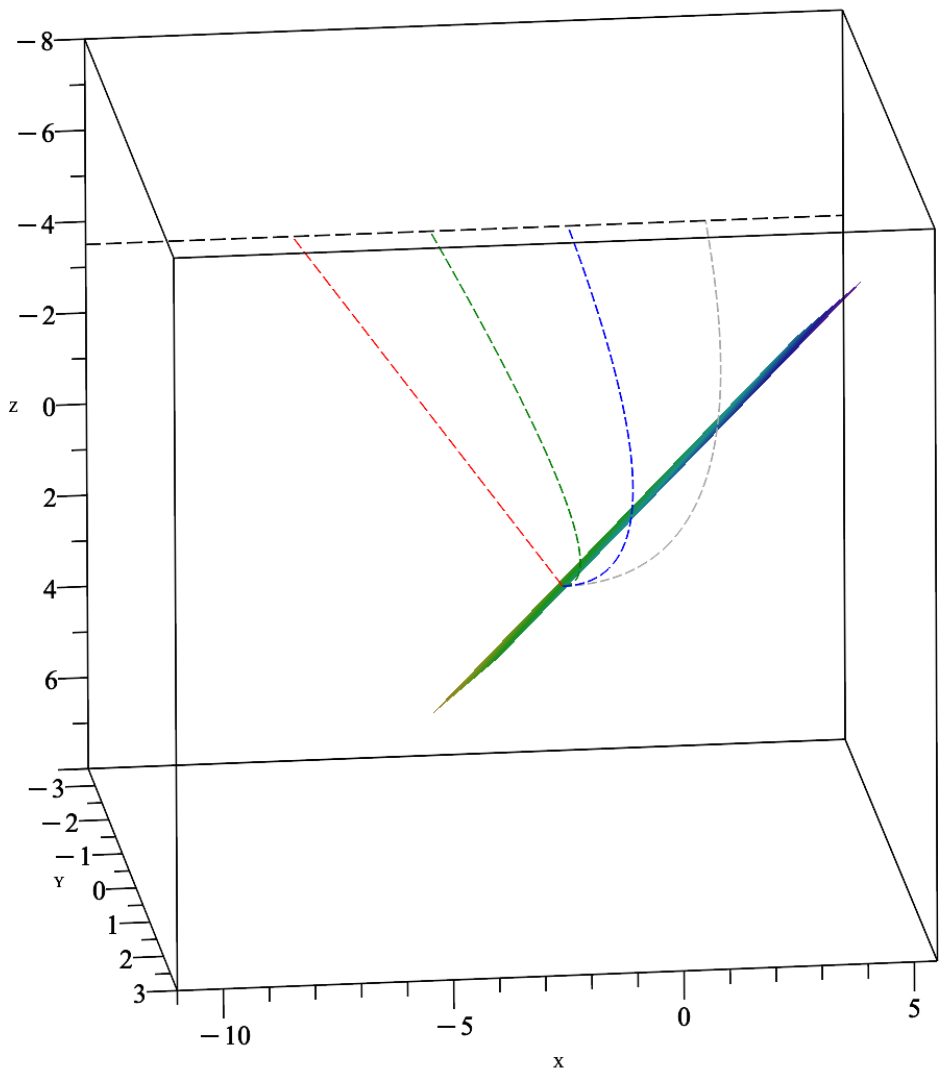}
	\caption{For fixed parameters $\rho_0=1$, $k_0=1$, $m_0=1$, we depict trajectories~\eqref{eqn7} of four particles starting  from the same initial point   $(-2,1,1)$ lying on the plane \eqref{eqn9}, and emitted with four different initial velocities $u_0 = 0, 1, 2, 3$. At $t=3$, they are shown to lie  on the same black dashed line parallel to the $x$-axis. }\label{fig2}
\end{figure}
The transformation $t \to -t$, $u_0 \to -u_0$ preserves the form of the world line equations for particle motion~\eqref{eqn7}. Consequently, it is sufficient to analyze the motion of particles for $t > 0$.
The Jacobian determinant of the transformation~\eqref{eqn7} from the Lagrangian to Euler coordinates is
\begin{equation*}
	\begin{array}{c}
		J=\begin{vmatrix}
			t& -k_0 & -m_0\\
			0& 1& 0\\
			0& 0& 1
		\end{vmatrix}
	\end{array}=t.
\end{equation*}
If the constants $k_0$ and $m_0$ are fixed, then from \eqref{eqn5_1} we see that, at $t=0$, all particles are located on the plane:
\begin{equation}\label{eqn9}
    x+k_0y + m_0z=0,
\end{equation}
the motion of the media is non-isochoric at any other moment of time. As time increases the volume of the particles increases as they  disperse throughout all space as $t\to\infty$.

If both $k_0=0$ and $m_0=0$, then the particles are moving along a line parallel to the $x$-axis. Otherwise, the projection of the particles onto the $yz$-plane is  a straight line:
\begin{equation*}
{m_0}(y-y_0)-k_0(z-z_0)=0.
\end{equation*}
If $k_0\neq 0$ ($m_0\neq 0$) then the projection of the trajectories onto the $xy$-plane ($xz$-plane) is a parabola.
 Particles emitted from  the same initial $(x_0,y_0, z_0)$, satisfying \eqref{eqn9}, with varying initial velocities $u_0$, at a fixed time $t$, are located on  the same line parallel to the $x$-axis.
Some trajectories of the particles are plotted in Figure\ref{fig2}.

\section{Conclusion}
In continuation of the study of four-dimensional subalgebras of the full 12-dimensional symmetry Lie algebra, admitted by the gas dynamics system with a special equation of state \eqref{state-sp}, we made the following progress. We computed invariants for a large set of subalgebras from an optimal list, determined the isomorphism classes of these subalgebras, and explicitly constructed a partially reduced system (partially invariant submodel) for one of these subalgebras. This led to a family of explicit solutions for the original system, whose trajectories we analyzed.  Our contribution provides a seed for future work on constructing other invariant and partially invariant submodels,  analyzing the hierarchy of submodels, and obtaining new explicit solutions.  
\bibliographystyle{plain}  
\bibliography{references}  

\end{document}